\documentclass[12pt]{article}
\usepackage{latexsym, amssymb, amsmath, amscd, amsfonts, epsfig, graphicx, colordvi,verbatim,ifpdf,extarrows}
\usepackage{amsfonts, amsmath, amssymb}
\usepackage{amssymb,amsfonts,amsmath,latexsym,epsfig,cite, psfrag,eepic,color}
\usepackage{amscd,graphics}
\usepackage{latexsym, amssymb,  amsmath,amscd, amsfonts, epsfig, graphicx, colordvi,amsthm}
\usepackage{amsmath}
\usepackage{graphicx}
\usepackage{epstopdf}
\usepackage{color}
\usepackage{ifpdf}
\usepackage{fancybox}
\usepackage[font=small,labelfont=bf,labelsep=none]{caption}
\usepackage{float}

\allowdisplaybreaks

\newtheorem{thm}{Theorem}[section]

\newtheorem{conj}[thm]{Conjecture}

\newtheorem{lem}[thm]{Lemma}
\newtheorem{core}[thm]{Corollary}

\def\pf{\noindent{\it Proof.} }
\def\qed{\nopagebreak\hfill{\rule{4pt}{7pt}}
\medbreak}

\numberwithin{equation}{section}

\def\qed{\nopagebreak\hfill{\rule{4pt}{7pt}}
\medbreak}

\newlength{\boxedparwidth}
  {\begin{center} \begin{tabular}{|@{\hspace{.315in}}c@{\hspace{.15in}}|}
                  \hline \\ \begin{minipage}[t]{\boxedparwidth}
                  \setlength{\parindent}{.25in}}%
  {\end{minipage} \\ \\ \hline \end{tabular} \end{center}}

\begin{document}
\begin{center}

 {\Large \bf Minimal excludant integer and bilateral truncated Jacobi triple product identity}
\end{center}

\begin{center}
{Y.H. Chen}$^{1}$, {W.D. Deng}$^{2}$, {Thomas Y. He}$^{3}$ and
  {H.X. Huang}$^{4}$ \vskip 2mm

$^{1,2,3,4}$ School of Mathematical Sciences, Sichuan Normal University, Chengdu 610066, P.R. China

   \vskip 2mm

$^1$chenyh@stu.sicnu.edu.cn, $^2$wddeng@stu.sicnu.edu.cn,  $^3$heyao@sicnu.edu.cn,  $^4$huanghaoxuan@stu.sicnu.edu.cn
\end{center}

\vskip 6mm   {\noindent \bf Abstract.} In 2012, Andrews and Merca proved a truncated theorem on Euler's pentagonal number theorem. Since then, a number of results on truncated theta series have been proved, including truncated Jacobi triple product identity. In this paper, we provide  partition-theoretic interpretations for the bilateral truncated Jacobi triple product identity in terms of the minimal excludant integer.


\section{Introduction}

A partition $\pi$ of a positive integer $n$ is a finite weakly decreasing sequence of positive integers $\pi=(\pi_1,\pi_2,\ldots,\pi_\ell)$ such that $\pi_1+\pi_2+\cdots+\pi_\ell=n$. The $\pi_i$ are called the parts of $\pi$. The empty sequence forms the only partition of zero. We use $|\pi|$ to denote the sum of the parts of $\pi$. Let $p(n)$ be the number of partitions of $n$. It is known that the generating function of $p(n)$  is
\begin{equation*}\label{p}
\sum_{n=0}^\infty p(n)q^n=\frac{1}{(q;q)_\infty}.
\end{equation*}
Here and in the sequel, we assume that $|q|<1$ and use the standard notation \cite{Andrews-1976}:
\[(a;q)_\infty=\prod_{i=0}^{\infty}(1-aq^i),\]
\[(a_1,a_2,\ldots,a_m;q)_\infty=(a_1;q)_\infty(a_2;q)_\infty\cdots(a_m;q)_\infty,\]
\[(a;q)_n=\frac{(a;q)_\infty}{(aq^n;q)_\infty},\]
and
\[{M\brack N}_{q^k}=\left\{\begin{array}{ll}\frac{(q^k;q^k)_M}{(q^k;q^k)_N(q^k;q^k)_{M-N}},&\text{if }M\geq N\geq 0,\\
0,&\text{otherwise.}
\end{array}\right.\]

The reciprocal of the generating function of $p(n)$ has the following series expansion:
\begin{equation}\label{pn}
(q;q)_\infty=\sum_{n=-\infty}^\infty(-1)^nq^{n(3n-1)/2},
\end{equation}
which is known as Euler's pentagonal number theorem.

In 2012, Andrews and Merca \cite{Andrews-Merca-2012} obtained a truncated version of  \eqref{pn}. For $k\geq 1$,
\begin{equation*}\label{pn1}
\frac{1}{(q;q)_\infty} \sum_{j=-k+1}^{k}(-1)^jq^{j(3j-1)/2}=1+(-1)^{k-1}\sum_{n=k}^\infty \frac{q^{k(k-1)/2+(k+1)n}}{(q;q)_n} \begin{bmatrix} n-1 \\ k-1 \end{bmatrix}_{q},
\end{equation*}
from which they deduced the following partition theorem.
\begin{thm}\cite[Theorem 1.1]{Andrews-Merca-2012}
For $n\geq 1$ and $k\ge1$,
\[
(-1)^{k-1}\sum_{j=-k+1}^{k}(-1)^jp(n-j(3j-1)/2)=M_k(n),
\]
where $M_k(n)$ is the number of partitions of $n$ where $k$ is the least integer that does not occur as a part and there are more parts greater than $k$ than there are less than $k$.
\end{thm}

After the work of Andrews and Merca, Guo and Zeng \cite{Guo-Zeng-2013} established truncated versions of two identities of Gauss \cite[(2.2.12),(2.2.13)]{Andrews-1976}.
Andrews and Merca \cite[Question (2)]{Andrews-Merca-2012} and Guo and Zeng\cite[Conjecture 6.1]{Guo-Zeng-2013}   considered truncated version of Jacobi triple product identity
\begin{equation}\label{jacobi-rs}
(q^S,q^{R-S},q^R;q^R)_\infty=\sum_{n=-\infty}^{\infty}(-1)^nq^{Rn(n-1)/2+Sn},
\end{equation}
and proposed the following conjecture.
\begin{conj}\label{amgz}
For integers $k$, $R$, and $S$ with $k\ge1$ and $1\le S\le R/2$, the coefficient of $q^n$ with $n\ge 1$ in
\[
\frac{(-1)^{k-1}}{(q^S,q^{R-S},q^R;q^R)_\infty}\sum_{j=-k+1}^{k}(-1)^jq^{Rj(j-1)/2+Sj}
\]
is nonnegative.
\end{conj}

Conjecture \ref{amgz} was proved independently by Mao \cite{Mao-2015} and Yee \cite{Yee-2015}, and then reconfirmed by Wang and Yee \cite{Wang-Yee-2019} by providing an explicit series form with nonnegative coefficients. Wang and Yee \cite[Theorems 1.1 and 1.2]{Wang-Yee-2019} obtained that for $1\le S\le R/2$ and $m\ge 1$,
\begin{equation}\label{jacobi-eqn-1}
\begin{split}
&\quad\frac{1}{(q^S,q^{R-S},q^R;q^R)_\infty}\sum_{n=-m+1}^{m}(-1)^nq^{Rn(n-1)/2+Sn}\\
&=1+(-1)^{m-1}q^{\binom{m}{2}R}\sum_{n=m}^{\infty}\sum_{\substack{i+j+h+k=n\\i,j,h,k\ge 0}}\frac{q^{(mj+hk)R+(h-k)S+nR}}{(q^R;q^R)_i(q^R;q^R)_j(q^R;q^R)_h(q^R;q^R)_k}
\left[\begin{array}{c}n-1\\m-1\end{array}\right]_{q^R},
\end{split}
\end{equation}
and
\begin{equation}\label{jacobi-eqn-2}
\begin{split}
&\quad\frac{1}{(q^S,q^{R-S},q^R;q^R)_\infty}\sum_{n=-m}^{m}(-1)^nq^{Rn(n-1)/2+Sn}\\
&=1+(-1)^m q^{\binom{m+1}{2}R}\sum_{n=m+1}^{\infty}\sum_{\substack{i+j+h+k=n\\i,j,h,k\ge 0}}\frac{q^{(mj+h(k-1))R+(h-k)S+nR}}{(q^R; q^R)_i(q^R;q^R)_j(q^R;q^R)_h(q^R;q^R)_k}
\left[\begin{array}{c}n-1\\ m \end{array}\right]_{q^R}.
\end{split}
\end{equation}

For $1\leq S<R$ and $n\geq 0$, let $J_{R,S}(n)$ be the number of triples $(\alpha,\beta,\gamma)$ of partitions such that $|\alpha|+|\beta|+|\gamma|=n$ and the parts of $\alpha$, $\beta$ and $\gamma$ are congruent to $S$, $-S$ and $0$ modulo $R$, respectively.  Clearly, the generating function of $J_{R,S}(n)$ is
\begin{equation*}\label{jacobi}
\sum_{n=0}^{\infty}J_{R,S}(n)q^{n}=\frac{1}{(q^{S},q^{R-S},q^{R};q^{R})_{\infty}}.
\end{equation*}
It follows from \eqref{jacobi-eqn-1} and \eqref{jacobi-eqn-2} that for $1\le S\le R/2$, $m\ge 1$ and $n\geq 1$,
\begin{equation}\label{li_1}
(-1)^{m-1}\sum_{j=-m+1}^{m}(-1)^j J_{R,S}\left(n-Rn(n-1)/2-Sn\right)\geq0,
\end{equation}
and
\begin{equation}\label{li_2}
(-1)^m\sum_{j=-m}^m(-1)^j J_{R,S}\left(n-Rn(n-1)/2-Sn\right)\ge 0.
\end{equation}

In 2021, Merca\cite{Merca-2021} considered another new truncated version of Euler's pentagonal number theorem and deduced that for $n\ge 1$ and $k\ge0$,
\begin{equation}\label{m1}
(-1)^k\sum_{j=-k}^k(-1)^j p(n-j(3j-1)/2)\ge 0.
\end{equation}

Xia and Zhao\cite{Xia-Zhao-2023} obtained the following partition-theoretic interpretation for the left-hand side of (\ref{m1}).
For $n\ge 1$ and $k\ge 0$,
\begin{equation*}\label{bmkn}
(-1)^k \sum_{j=-k}^k(-1)^jp(n-j(3j-1)/2)=\widetilde{P}_k(n),
\end{equation*}
where $\widetilde{P}_k(n)$ is the number of partitions of $n$ in which $k$ is the least integer such that every
part $\le k$ appears at least once and the first part larger than $k$ appears at least $k+1$ times.

Furthermore, Xia and Zhao\cite{Xia-Zhao-2023} proved that for $m\le k$ and $n\geq 1$,
\begin{equation*}
(-1)^{min\{|m|,k\}}\sum_{j=m}^{k}(-1)^jp(n-j(3j-1)/2)\ge 0.
\end{equation*}

Then, Li \cite{Li-2023} generalized the truncated sums of Jacobi triple product series in (\ref{li_1}) and (\ref{li_2}).
\begin{thm}\cite[Corollary 5.1]{Li-2023}\label{lit1}
For integers \(R\), \(S\), \(m\), \(k\), and \(n\) with
\(1\leq S\leq R/2\), \(m\leq k\), and \(n\geq 1\),
\begin{equation*}
(-1)^{min\{|m|,k\}}\sum_{j=m}^{k}(-1)^{j}J_{R,S}(n-Rj(j-1)/2-Sj)\geq 0.
\end{equation*}
\end{thm}

Merca \cite{Merca-2021a} investigated truncated Jacobi triple product identity written in the form
\begin{equation}\label{stronger-jacobi-000}
\frac{1}{(q^S,q^{R-S};q^R)_\infty}\sum_{n=-\infty}^{\infty}(-1)^nq^{Rn(n-1)/2+Sn}-(q^R;q^R)_\infty=0,
\end{equation}
and made the following conjecture.

\begin{conj}\label{stronger-jacobi-con}
For integers $k$, $R$, and $S$ with $k\ge1$ and $1\le S<R$, the coefficient of $q^n$ with $n\ge 1$ in
\[
\frac{(-1)^{k}}{(q^S,q^{R-S};q^R)_\infty}\sum_{j=k}^{\infty}(-1)^jq^{Rj(j+1)/2-Sj}(1-q^{S(2j+1)})
\]
is nonnegative.
\end{conj}

Ballantine and Feigon \cite{Ballantine-Feigon-2012} proved Conjecture \ref{stronger-jacobi-con} for $k\in\{1, 2, 3\}$.
Ding and Sun showed that Conjecture \ref{stronger-jacobi-con} holds for $R = 2S$ in \cite{Ding-Sun-2025a} and $R = 3S$ in \cite{Ding-Sun-2025}, and
provided a systematic method
to determine an integer $N(R,S,k)$ such that Conjecture \ref{stronger-jacobi-con} holds for $n\geq N(R,S,k)$ in \cite{Ding-Sun-2025a}.

Conjecture \ref{stronger-jacobi-con} is equivalent to showing that for $k\ge1$ and $1\le S<R$, the coefficient of $q^n$ with $n\ge 1$ in
\begin{equation*}
\frac{(-1)^{k-1}}{(q^S,q^{R-S};q^R)_\infty}\sum_{j=-k+1}^{k}(-1)^jq^{Rj(j-1)/2+Sj}-(-1)^{k-1}(q^R;q^R)_\infty
\end{equation*}
is nonnegative.

For integers $R$ and $S$ with $1\leq S<R$, and integers $m$ and $k$ with $m\leq k$, we consider the following bilateral truncated version of \eqref{stronger-jacobi-000}:
\begin{equation}\label{eqn-main-result}
(-1)^{min\{|m|,k\}}\left(\frac{1}{(q^S,q^{R-S};q^R)_\infty}\sum_{j=m}^{k}(-1)^jq^{Rj(j-1)/2+Sj}-(q^R;q^R)_\infty\right).
\end{equation}
The main objective of this article is to provide partition-theoretic interpretations of \eqref{eqn-main-result} by using the minimal excludant integer.

Recently, Andrews and Newman \cite{Andrews-2019} undertook a combinatorial study of the minimal excludant of a partition, which was earlier introduced by Grabner and Knopfmacher \cite{Grabner-Knopfmacher-2006} under the name ``smallest gap". The minimal excludant of a partition $\pi$ is the smallest positive integer that is not a part of $\pi$. For example, among the partitions of $4$, the minimal excludants of $(4)$, $(3,1)$, $(2,2)$, $(2,1,1)$, and $(1,1,1,1)$ are $1$, $2$, $1$, $3$, and $2$, respectively.

In \cite{Andrews-Newman-2020},  Andrews and Newman gave the definition of the minimal excludant of a partition in congruence classes. For $A\geq a\geq 1$,  they defined $mex_{A,a}(\pi)$ to be the smallest positive integer congruent to $a$ modulo $A$ that is not a part of $\pi$. For example, consider $A=2$ and $a=1$. Then $mex_{2,1}(\pi)$ is the smallest odd positive integer that is not a part of $\pi$. We have $mex_{2,1}((4))=1$, $mex_{2,1}((3,1))=5$, $mex_{2,1}((2,2))=1$, $mex_{2,1}((2,1,1))=3$, and $mex_{2,1}((1,1,1,1))=3$.

For $1\leq S< R$ and $n\geq 0$, let $\mathcal{P}_{R,S}(n)$ be the set of pairs $(\alpha,\beta)$ of partitions such that $|\alpha|+|\beta|=n$ and the parts of $\alpha$ and $\beta$ are congruent to $S$ and $-S$ modulo $R$, respectively.  Let ${P}_{R,S}(n)$ be the number of partitions in $\mathcal{P}_{R,S}(n)$, Clearly, the generating function for the partitions in  $\mathcal{P}_{R,S}(n)$ is
\begin{equation}\label{jacobi00000000}
\sum_{n=0}^{\infty}P_{R,S}(n)q^{n}=\frac{1}{(q^{S},q^{R-S};q^{R})_{\infty}}.
\end{equation}

For $d\geq 0$, we use ${P}_{R,S,d}(n)$ (resp. ${P}_{R,-S,d+1}(n)$) to denote the number of pairs $(\alpha,\beta)$ in $\mathcal{P}_{R,S}(n)$ such that
$mex_{R,S}(\alpha)\geq Rd+S$ (resp. $mex_{R,R-S}(\beta)\geq R(d+1)-S$) and $mex_{R,S}(\alpha)\equiv Rd+S\pmod{2R}$ (resp. $mex_{R,R-S}(\beta)\equiv R(d+1)-S\pmod{2R}$).

Assume that
\begin{equation}\label{gen-rho}
\sum_{n=0}^\infty \rho_R(n)q^n=(q^R;q^R)_\infty,
\end{equation}
we give the following partition-theoretic interpretations of \eqref{eqn-main-result}.
\begin{thm}\label{main-thm} Let $R,S$ and $n$ be integers such that $1\leq S<R$ and $n\geq 0$.
For $m<0<k$, we have
\begin{equation}\label{mk1}
\sum_{j=m}^{k}(-1)^jP_{R,S}(n-Rj(j-1)/2-Sj)-\rho_R(n)=(-1)^k{P}_{R,S,k+1}(n)+(-1)^m{P}_{R,-S,2-m}(n).
\end{equation}
For $k\geq m\geq 0$, we have
\begin{equation}\label{mk2}
(-1)^{min\{|m|,k\}}\sum_{j=m}^{k}(-1)^jP_{R,S}(n-Rj(j-1)/2-Sj)={P}_{R,S,m}(n)+(-1)^{m+k}{P}_{R,S,k+1}(n).
\end{equation}
For $m\leq k\leq 0$, we have
\begin{equation}\label{mk3}
(-1)^{min\{|m|,k\}}\sum_{j=m}^{k}(-1)^jP_{R,S}(n-Rj(j-1)/2-Sj)={P}_{R,-S,1-k}(n)+(-1)^{m+k}{P}_{R,-S,2-m}(n).
\end{equation}
\end{thm}

As a corollary of Theorem \ref{main-thm}, we can get the following result.
\begin{core}\label{main-core-proof}
Let $R$, $S$ and $n$ be integers such that $1\le S<R$ and $n\geq 0$. For $m<0<k$ and $m\equiv k\pmod{2}$, we have
\begin{equation}\label{proof-mk-0}
(-1)^{min\{|m|,k\}}\left(\sum_{j=m}^{k}(-1)^jP_{R,S}(n-Rj(j-1)/2-Sj)-\rho_R(n)\right)\geq 0.
\end{equation}
For  $m\le k$ and $mk\geq 0$, we have
\begin{equation}\label{proof-mk-1-3}
(-1)^{min\{|m|,k\}}\sum_{j=m}^{k}(-1)^jP_{R,S}(n-Rj(j-1)/2-Sj)\geq 0.
\end{equation}
\end{core}

For $1\leq S<R$ and $n\geq 0$, let $\mathcal{J}_{R,S}(n)$  be the set of triples $(\alpha,\beta,\gamma)$ of partitions such that $|\alpha|+|\beta|+|\gamma|=n$ and the parts of $\alpha$, $\beta$ and $\gamma$ are congruent to $S$, $-S$ and $0$ modulo $R$, respectively. Then,  $J_{R,S}(n)$ is the number of triples  of partitions in  $\mathcal{J}_{R,S}(n)$. For $d\geq 0$, we use ${J}_{R,S,d}(n)$ (resp. ${J}_{R,-S,d+1}(n)$) to denote the number of triples $(\alpha,\beta,\gamma)$ in $\mathcal{J}_{R,S}(n)$ such that
$mex_{R,S}(\alpha)\geq Rd+S$ (resp. $mex_{R,R-S}(\beta)\geq R(d+1)-S$) and $mex_{R,S}(\alpha)\equiv Rd+S\pmod{2R}$ (resp. $mex_{R,R-S}(\beta)\equiv R(d+1)-S\pmod{2R}$). Note that
\begin{equation}\label{add-label*}\sum_{n=0}^\infty J_{R,S}(n)q^n=\frac{1}{(q^R;q^R)_\infty}\sum_{n=0}^\infty P_{R,S}(n)q^n,\end{equation}
\[\sum_{n=0}^\infty {J}_{R,S,d}(n)q^n=\frac{1}{(q^R;q^R)_\infty}\sum_{n=0}^\infty {P}_{R,S,d}(n)q^n,\]
and
\[\sum_{n=0}^\infty {J}_{R,-S,d+1}(n)q^n=\frac{1}{(q^R;q^R)_\infty}\sum_{n=0}^\infty {P}_{R,-S,d+1}(n)q^n,\]
we can get the following corollaries of Theorem \ref{main-thm}.

\begin{core}\label{add*core-01} Let $R$ and $S$  be integers with $1\leq S<R$.
For $m<0<k$ and $n\geq 1$, we have
\begin{equation*}\label{mkJ1}
\sum_{j=m}^{k}(-1)^jJ_{R,S}(n-Rj(j-1)/2-Sj)=(-1)^k{J}_{R,S,k+1}(n)+(-1)^m{J}_{R,-S,2-m}(n).
\end{equation*}
For $k\geq m\geq 0$ and $n\geq 0$, we have
\begin{equation*}\label{mkJ2}
(-1)^{min\{|m|,k\}}\sum_{j=m}^{k}(-1)^jJ_{R,S}(n-Rj(j-1)/2-Sj)={J}_{R,S,m}(n)+(-1)^{m+k}{J}_{R,S,k+1}(n).
\end{equation*}
For $m\leq k\leq 0$ and $n\geq 0$, we have
\begin{equation*}\label{mkJ3}
(-1)^{min\{|m|,k\}}\sum_{j=m}^{k}(-1)^jJ_{R,S}(n-Rj(j-1)/2-Sj)={J}_{R,-S,1-k}(n)+(-1)^{m+k}{J}_{R,-S,2-m}(n).
\end{equation*}
\end{core}

\begin{core}\label{add*core}
For integers $R$ and $S$ with $1\le S<R$, and integers $m$, $k$, and $n$ satisfying either $m<0<k$, $m\equiv k\pmod{2}$  and $n\geq 1$, or $m\le k$, $mk\geq 0$ and $n\geq 0$,
\begin{equation*}
(-1)^{min\{|m|,k\}}\sum_{j=m}^{k}(-1)^jJ_{R,S}(n-Rj(j-1)/2-Sj)\geq 0.
\end{equation*}
\end{core}

For integers $R$, $S$ and $n$ with $1\le S<R$ and $n\geq 0$, and integers $m$ and $k$ satisfying either $m<0<k$ and $m\equiv k\pmod{2}$, or $m\le k$ and $mk\geq 0$ , we set
\[(-1)^{min\{|m|,k\}}\sum_{j=m}^{k}(-1)^jJ_{R,S}(n-Rj(j-1)/2-Sj)=C_{m,k}(n).\]
It follows from Corollary \ref{add*core} that $C_{m,k}(n)\geq 0$ when $m<0<k$, $m\equiv k\pmod{2}$  and $n\geq 1$, or $m\le k$, $mk\geq 0$ and $n\geq 0$. Then, we have the following result.
\begin{thm}\label{main-thm-proof-new} Let $R$, $S$, $\ell$ and $n$ be integers with $1\le S<R$, $\ell\geq 1$ and $n\geq 0$. For $m<0<k$ and $m\equiv k\pmod{2}$, we have
\begin{equation*}
\sum_{i=-\infty}^{\infty}(-1)^iC_{m,k}(n-\ell Ri(3i-1)/2)-(-1)^k\rho_{R\ell}(n)\geq 0.
\end{equation*}
For $m\le k$ and $mk\geq 0$, we have
\begin{equation*}
\sum_{i=-\infty}^{\infty}(-1)^iC_{m,k}(n-\ell Ri(3i-1)/2)\geq 0.
\end{equation*}
\end{thm}

By Theorem \ref{lit1} with $m=-k+1<0$ and Corollary \ref{add*core-01} with $m=-k+1<0$, we get the following result.

\begin{core}
Let \(R\), \(S\), \(k\), and \(n\) be integers  with
\(1\leq S\leq R/2\), \(k\geq 2\), and \(n\geq 1\). Then
\begin{equation*}
J_{R,-S,k+1}(n)\geq J_{R,S,k+1}(n).
\end{equation*}
\end{core}

By definition,  we can get that for $k\geq 2$,
\begin{equation}\label{sum-side}
 \begin{split}
 &\sum_{n=0}^\infty \left(J_{R,-S,k+1}(n)-J_{R,S,k+1}(n)\right)q^n\\
 &\quad=\frac{1}{(q^S,q^{R-S},q^R; q^R)_\infty}\left(\sum_{j=0}^{\infty} q^{R(2j+k+1)(2j+k)/2 - S(2j+k)} (1 - q^{R(2j+k+1)-S})\right.\\
&\qquad\qquad\qquad\qquad\qquad\quad\left.-\sum_{j=0}^{\infty} q^{R(2j+k+1)(2j+k)/2+S(2j+k+1)} (1 - q^{R(2j+k+1)+S}) \right).
\end{split}
\end{equation}

We now focus on the summand appearing on the right-hand side of \eqref{sum-side}, namely
\[q^{R(2j+k+1)(2j+k)/2 - S(2j+k)} (1 - q^{R(2j+k+1)-S})- q^{R(2j+k+1)(2j+k)/2 + S(2j+k+1)} (1 - q^{R(2j+k+1)+S}).\]
For simplicity, let $t=2j+k$. We then study
\[q^{R(t+1)t/2 - St} (1 - q^{R(t+1)-S})-q^{R(t+1)t/2 + S(t+1)} (1 - q^{R(t+1)+S}).\]

\begin{thm}\label{main-thm-proof-FHG}
For $t\geq 2$, we have
\begin{equation*}
\begin{split}
&\quad q^{R(t+1)t/2 - St} (1 - q^{R(t+1)-S})-q^{R(t+1)t/2 + S(t+1)} (1 - q^{R(t+1)+S})\\
&=q^{R(t+1)t/2}\sum_{i=0}^{t-1} H_i+q^{R(t+1)t/2+R(t+1)t}(1-q^S)-q^{R(t+1)t/2+R(t+1)-S(t+1)}(1-q^{S(2t+3)}),
\end{split}
\end{equation*}
where
\begin{align*}
H_i&=q^{iR(t+1)-S(t-i)}\left(1-q^{2S}+(q^{2S}-q^{S(2t+1-2i)})(1-q^{R(t+1)-S})\right)\\
&=q^{iR(t+1)-S(t-i)}\left(1-q^{R(t+1)+S}-q^{S(2t+1-2i)}+q^{R(t+1)+S(2t-2i)}\right).
\end{align*}
\end{thm}

This article is organized as follows. We will show Theorem \ref{main-thm}, Corollary \ref{main-core-proof} and Theorem \ref{main-thm-proof-new} in Section 2 and prove Theorem \ref{main-thm-proof-FHG} in Section 3.

\section{Proofs of Theorem \ref{main-thm}, Corollary \ref{main-core-proof} and Theorem \ref{main-thm-proof-new}}

In this section, we first prove Theorem \ref{main-thm} and then derive Corollary \ref{main-core-proof} from it. Finally,  with the aid of Corollary \ref{main-core-proof}, we prove Theorem \ref{main-thm-proof-new}.

Before proving Theorem \ref{main-thm}, we are required to give the generating function of ${P}_{R,S,d}(n)$ and ${P}_{R,-S,d+1}(n)$.
\begin{lem}\label{lem-main-proof}
For $1\leq S< R$ and $d\geq 0$, we have
\begin{equation}\label{lem-eqn-1}
\sum_{n=0}^\infty {P}_{R,S,d}(n)q^n=\frac{1}{(q^S,q^{R-S};q^R)_\infty}\sum_{j=0}^{\infty}q^{R(2j+d)(2j+d-1)/2+S(2j+d)}(1-q^{R(2j+d)+S})
\end{equation}
and
\begin{equation}\label{lem-eqn-2}
\sum_{n=0}^\infty {P}_{R,-S,d+1}(n)q^n=\frac{1}{(q^S,q^{R-S};q^R)_\infty}\sum_{j=0}^{\infty}q^{R(2j+d+1)(2j+d)/2-S(2j+d)}(1-q^{R(2j+d+1)-S}).
\end{equation}
\end{lem}

\pf For $n\geq 0$, suppose that  $(\alpha,\beta)$ is a pair counted by ${P}_{R,S,d}(n)$, then by definition, we have $mex_{R,S}(\alpha)=R(2j+d)+S$ for some integer $j\geq 0$. By the definition of \(mex_{R,S}(\alpha)\), we see that the parts
\(S,R+S,\ldots,R(2j+d-1)+S\) appear in \(\alpha\), whereas
\(R(2j+d)+S\) does not occur in \(\alpha\).  Therefore, the generating function of ${P}_{R,S,d}(n)$ is
\begin{equation*}
\begin{split}
\sum_{n=0}^\infty {P}_{R,S,d}(n)q^n&=\frac{1}{(q^{R-S};q^R)_\infty}\sum_{j=0}^{\infty}\frac{q^{S+(R+S)+\cdots+(R(2j+d-1)+S)}(1-q^{R(2j+d)+S})}{(q^S;q^R)_\infty}\\
&=\frac{1}{(q^S,q^{R-S};q^R)_\infty}\sum_{j=0}^{\infty}q^{R(2j+d)(2j+d-1)/2+S(2j+d)}(1-q^{R(2j+d)+S}),
\end{split}
\end{equation*}
and thus \eqref{lem-eqn-1} is valid.  With a similar argument above, we can get

\begin{equation*}
\begin{split}
\sum_{n=0}^\infty {P}_{R,-S,d+1}(n)q^n&=\frac{1}{(q^{S};q^R)_\infty}\sum_{j=0}^{\infty}\frac{q^{R-S+(2R-S)+\cdots+(R(2j+d)-S)}(1-q^{R(2j+d+1)-S})}{(q^{R-S};q^R)_\infty}\\
&=\frac{1}{(q^S,q^{R-S};q^R)_\infty}\sum_{j=0}^{\infty}q^{R(2j+d+1)(2j+d)/2-S(2j+d)}(1-q^{R(2j+d+1)-S}).
\end{split}
\end{equation*}
We arrive at \eqref{lem-eqn-2}. This completes the proof.  \qed

Now, we proceed to prove Theorem \ref{main-thm}.

{\noindent \bf Proof of Theorem \ref{main-thm}.}
For $m<0<k$, we have
\begin{align*}
&\quad\sum_{n=0}^\infty\left(\sum_{j=m}^{k}(-1)^jP_{R,S}(n-Rj(j-1)/2-Sj)-\rho_R(n)\right)q^n\\
&=\sum_{j=m}^{k}(-1)^jq^{Rj(j-1)/2+Sj}\sum_{n=0}^\infty P_{R,S}(n)q^n-\sum_{n=0}^\infty\rho_R(n)q^n\\
&=\frac{1}{(q^S,q^{R-S};q^R)_\infty}\sum_{j=m}^{k}(-1)^jq^{Rj(j-1)/2+Sj}-(q^R;q^R)_\infty\\
&=\frac{1}{(q^S,q^{R-S};q^R)_\infty}\left(\sum_{j=-\infty}^{\infty}(-1)^jq^{Rj(j-1)/2+Sj}-\sum_{j=k+1}^{\infty}(-1)^jq^{Rj(j-1)/2+Sj}\right.\\
&\qquad\qquad\qquad\qquad\left.-\sum_{j=-\infty}^{m-1}(-1)^jq^{Rj(j-1)/2+Sj}\right)-(q^R;q^R)_\infty\\
&=\frac{1}{(q^S,q^{R-S};q^R)_\infty}\left((q^S,q^{R-S},q^R;q^R)_\infty-\sum_{j=0}^{\infty}(-1)^{j+k+1}q^{R(j+k+1)(j+k)/2+S(j+k+1)}\right.\\
&\qquad\qquad\qquad\qquad\left.-\sum_{j=0}^{\infty}(-1)^{-j+m-1}q^{R(-j+m-1)(-j+m-2)/2+S(-j+m-1)}\right)-(q^R;q^R)_\infty\\
&=\frac{(-1)^k}{(q^S,q^{R-S};q^R)_\infty}\sum_{j=0}^{\infty}(-1)^{j}q^{R(j+k+1)(j+k)/2+S(j+k+1)}\\
&\quad+\frac{(-1)^m}{(q^S,q^{R-S};q^R)_\infty}\sum_{j=0}^{\infty}(-1)^{j}q^{R(j-m+1)(j-m+2)/2-S(j-m+1)}\\
&=\frac{(-1)^k}{(q^S,q^{R-S};q^R)_\infty}\left(\sum_{j=0}^{\infty}q^{R(2j+k+1)(2j+k)/2+S(2j+k+1)}\right.\\
&\qquad\qquad\qquad\qquad\quad\left.-\sum_{j=0}^{\infty}q^{R(2j+k+2)(2j+k+1)/2+S(2j+k+2)}\right)\\
&\quad+\frac{(-1)^m}{(q^S,q^{R-S};q^R)_\infty}\left(\sum_{j=0}^{\infty}q^{R(2j-m+1)(2j-m+2)/2-S(2j-m+1)}\right.\\
&\qquad\qquad\qquad\qquad\qquad\left.-\sum_{j=0}^{\infty}q^{R(2j-m+2)(2j-m+3)/2-S(2j-m+2)}\right)\\
&=\frac{(-1)^k}{(q^S,q^{R-S};q^R)_\infty}\sum_{j=0}^{\infty}q^{R(2j+k+1)(2j+k)/2+S(2j+k+1)}(1-q^{R(2j+k+1)+S})\\
&\quad+\frac{(-1)^m}{(q^S,q^{R-S};q^R)_\infty}\sum_{j=0}^{\infty}q^{R(2j-m+1)(2j-m+2)/2-S(2j-m+1)}(1-q^{R(2j-m+2)-S})\\
&=(-1)^k\sum_{n=0}^\infty{P}_{R,S,k+1}(n)q^n+(-1)^m\sum_{n=0}^\infty{P}_{R,-S,2-m}(n)q^n,
\end{align*}
where the second equality follows from \eqref{jacobi00000000} and \eqref{gen-rho}, the fourth equality follows from \eqref{jacobi-rs}, and
the final equality follows from Lemma \ref{lem-main-proof}. Hence  \eqref{mk1} follows.

 We then prove \eqref{mk2} and \eqref{mk3} by using a similar argument above.
For $k\geq m\geq 0$, in such case, we have $min\{|m|,k\}=m$. Using \eqref{jacobi00000000} and \eqref{lem-eqn-1}, we get
\begin{align*}
&\quad\sum_{n=0}^\infty\left((-1)^{min\{|m|,k\}}\sum_{j=m}^{k}(-1)^jP_{R,S}(n-Rj(j-1)/2-Sj)\right)q^n\\
&=(-1)^m\sum_{j=m}^{k}(-1)^jq^{Rj(j-1)/2+Sj}\sum_{n=0}^\infty P_{R,S}(n)q^n\\
&=\frac{(-1)^m}{(q^S,q^{R-S};q^R)_\infty}\left(\sum_{j=m}^{\infty}(-1)^jq^{Rj(j-1)/2+Sj}-\sum_{j=k+1}^{\infty}(-1)^jq^{Rj(j-1)/2+Sj}\right)\\
&=\frac{(-1)^m}{(q^S,q^{R-S};q^R)_\infty}\\
&\quad\times\left(\sum_{j=0}^{\infty}(-1)^{j+m}q^{R(j+m)(j+m-1)/2+S(j+m)}
-\sum_{j=0}^{\infty}(-1)^{j+k+1}q^{R(j+k+1)(j+k)/2+S(j+k+1)}\right)\\
&=\frac{1}{(q^S,q^{R-S};q^R)_\infty}\sum_{j=0}^{\infty}q^{R(2j+m)(2j+m-1)/2+S(2j+m)}(1-q^{R(2j+m)+S})\\
&\quad+\frac{(-1)^{m+k}}{(q^S,q^{R-S};q^R)_\infty}\sum_{j=0}^{\infty}q^{R(2j+k+1)(2j+k)/2+S(2j+k+1)}(1-q^{R(2j+k+1)+S})\\
&=\sum_{n=0}^\infty{P}_{R,S,m}(n)q^n+(-1)^{m+k}\sum_{n=0}^\infty{P}_{R,S,k+1}(n)q^n,
\end{align*}
which leads to \eqref{mk2}.

For $m\leq k\leq 0$, in this case, we have $min\{|m|,k\}=k$. It follows from \eqref{jacobi00000000} and \eqref{lem-eqn-2} that
\begin{align*}
&\quad\sum_{n=0}^\infty\left((-1)^{min\{|m|,k\}}\sum_{j=m}^{k}(-1)^jP_{R,S}(n-Rj(j-1)/2-Sj)\right)q^n\\
&=(-1)^k\sum_{j=m}^{k}(-1)^jq^{Rj(j-1)/2+Sj}\sum_{n=0}^\infty P_{R,S}(n)q^n\\
&=\frac{(-1)^k}{(q^S,q^{R-S};q^R)_\infty}\left(\sum_{j=-\infty}^{k}(-1)^jq^{Rj(j-1)/2+Sj}-\sum_{j=-\infty}^{m-1}(-1)^jq^{Rj(j-1)/2+Sj}\right)\\
&=\frac{(-1)^k}{(q^S,q^{R-S};q^R)_\infty}\left(\sum_{j=0}^{\infty}(-1)^{-j+k}q^{R(-j+k)(-j+k-1)/2+S(-j+k)}\right.\\
&\qquad\qquad\qquad\qquad\left.-\sum_{j=0}^{\infty}(-1)^{-j+m-1}q^{R(-j+m-1)(-j+m-2)/2+S(-j+m-1)}
\right)\\
&=\frac{1}{(q^S,q^{R-S};q^R)_\infty}\\
&\quad\times\left(\sum_{j=0}^{\infty}(-1)^{j}q^{R(j-k)(j-k+1)/2-S(j-k)}+(-1)^{m+k}\sum_{j=0}^{\infty}(-1)^{j}q^{R(j-m+1)(j-m+2)/2-S(j-m+1)}
\right)\\
&=\frac{1}{(q^S,q^{R-S};q^R)_\infty}\sum_{j=0}^{\infty}q^{R(2j-k)(2j-k+1)/2-S(2j-k)}(1-q^{R(2j-k+1)-S})\\
&\quad+\frac{(-1)^{m+k}}{(q^S,q^{R-S};q^R)_\infty}\sum_{j=0}^{\infty}q^{R(j-m+1)(j-m+2)/2-S(2j-m+1)}(1-q^{R(2j-m+2)-S})\\
&=\sum_{n=0}^\infty{P}_{R,-S,1-k}(n)q^n+(-1)^{m+k}\sum_{n=0}^\infty{P}_{R,-S,2-m}(n)q^n,
\end{align*}
and so \eqref{mk3} is valid. The proof is complete.  \qed

Next, we derive  Corollary \ref{main-core-proof} from Theorem \ref{main-thm}.

{\noindent\bf Proof of Corollary \ref{main-core-proof}.} Clearly, it follows from \eqref{mk1} that \eqref{proof-mk-0} holds for $m<0<k$ and $m\equiv k\pmod{2}$. Furthermore, by \eqref{mk2} and \eqref{mk3}, we see that \eqref{proof-mk-1-3} holds for $m\le k$, $mk\geq 0$ and $m\equiv k\pmod{2}$.
Again by \eqref{mk2} and \eqref{mk3}, for  $m\not\equiv k\pmod{2}$, we obtain
\begin{align*}
&\quad(-1)^{min\{|m|,k\}}\sum_{j=m}^{k}(-1)^jP_{R,S}(n-Rj(j-1)/2-Sj)\\
&=\left\{\begin{array}{ll}
{P}_{R,S,m}(n)-{P}_{R,S,k+1}(n),&\text{if }k\geq m\geq 0,\\
{P}_{R,-S,1-k}(n)-{P}_{R,-S,2-m}(n),&\text{if }m\leq k\leq 0.
\end{array}\right.
\end{align*}

By definition, we know that for $k\geq m\geq 0$, ${P}_{R,S,m}(n)-{P}_{R,S,k+1}(n)$ is the number of pairs $(\alpha,\beta)$ in $\mathcal{P}_{R,S}(n)$ such that $mex_{R,S}(\alpha)$ belongs to
\[\{Rm+S,R(m+2)+S,\ldots,R(k-1)+S\},\]
and for $m\leq k\leq 0$, ${P}_{R,-S,1-k}(n)-{P}_{R,-S,2-m}(n)$ is the number of pairs $(\alpha,\beta)$ in $\mathcal{P}_{R,S}(n)$ such that $mex_{R,S}(\beta)$ belongs to \[\{R(1-k)-S,R(3-k)-S,\ldots,R(-m)-S\}.\]
This implies that  \eqref{proof-mk-1-3} is valid for $m\le k$, $mk\geq 0$, and $m\not\equiv k\pmod{2}$. The proof is complete.  \qed

We conclude this section with a proof of Theorem \ref{main-thm-proof-new}.

{\noindent \bf Proof of Theorem \ref{main-thm-proof-new}.} The generating function of
\[\sum_{i=-\infty}^{\infty}(-1)^iC_{m,k}(n-\ell Ri(3i-1)/2)\]
is given by
\begin{align}
&\quad\sum_{n=0}^\infty\left(\sum_{i=-\infty}^{\infty}(-1)^iC_{m,k}(n-\ell Ri(3i-1)/2)\right)q^n\nonumber\\
&=\sum_{i=-\infty}^{\infty}(-1)^iq^{\ell Ri(3i-1)/2}\sum_{n=0}^\infty C_{m,k}(n)q^n\nonumber\\
&=(q^{\ell R};q^{\ell R})_\infty\sum_{n=0}^\infty \left((-1)^{min\{|m|,k\}}\sum_{j=m}^{k}(-1)^jJ_{R,S}(n-Rj(j-1)/2-Sj)\right)q^n\nonumber\\
&=\frac{(q^{\ell R};q^{\ell R})_\infty}{(q^R;q^R)_\infty}\sum_{n=0}^\infty \left((-1)^{min\{|m|,k\}}\sum_{j=m}^{k}(-1)^jP_{R,S}(n-Rj(j-1)/2-Sj)\right)q^n,\label{add*inter-proof-eqn}
\end{align}
where the second equality follows from \eqref{pn} with $q\rightarrow q^{\ell R}$ and the third equality follows from \eqref{add-label*}. Using \eqref{proof-mk-1-3}, we obtain that for $m\le k$ and $mk\geq 0$, the coefficient of $q^n$ with $n\geq 0$ in \eqref{add*inter-proof-eqn} is nonnegative.

In the case $m<0<k$ and $m\equiv k\pmod{2}$, we have 
\[(-1)^{min\{|m|,k\}}=(-1)^m=(-1)^k.\] 
Therefore, we get
\begin{align}
&\quad\sum_{n=0}^\infty \left((-1)^{min\{|m|,k\}}\sum_{j=m}^{k}(-1)^jP_{R,S}(n-Rj(j-1)/2-Sj)\right)q^n\nonumber\\
&=\sum_{n=0}^\infty (-1)^k\left(\sum_{j=m}^{k}(-1)^jP_{R,S}(n-Rj(j-1)/2-Sj)-\rho_R(n)+\rho_R(n)\right)q^n\nonumber\\
&=\sum_{n=0}^\infty (-1)^k\left(\sum_{j=m}^{k}(-1)^jP_{R,S}(n-Rj(j-1)/2-Sj)-\rho_R(n)\right)q^n+(-1)^k\sum_{n=0}^\infty\rho_R(n)q^n\nonumber\\
&=\sum_{n=0}^\infty (-1)^k\left(\sum_{j=m}^{k}(-1)^jP_{R,S}(n-Rj(j-1)/2-Sj)-\rho_R(n)\right)q^n+(-1)^k(q^R;q^R)_\infty,\label{add-label**}
\end{align}
where the final equality follows from \eqref{gen-rho}. Again by \eqref{gen-rho}, we get 
\begin{equation}\label{add-label***}
\frac{(q^{\ell R};q^{\ell R})_\infty}{(q^R;q^R)_\infty}(q^R;q^R)_\infty=(q^{\ell R};q^{\ell R})_\infty=\sum_{n=0}^\infty\rho_{\ell R}(n)q^n.
\end{equation}

Substituting \eqref{add-label**} into \eqref{add*inter-proof-eqn}, and using \eqref{add-label***}, we have

\begin{align}
&\quad\sum_{n=0}^\infty\left(\sum_{i=-\infty}^{\infty}(-1)^iC_{m,k}(n-\ell Ri(3i-1)/2)\right)q^n-(-1)^k\sum_{n=0}^\infty\rho_{\ell R}(n)q^n\nonumber\\
&=\frac{(q^{\ell R};q^{\ell R})_\infty}{(q^R;q^R)_\infty}\sum_{n=0}^\infty (-1)^k\left(\sum_{j=m}^{k}(-1)^jP_{R,S}(n-Rj(j-1)/2-Sj)-\rho_R(n)\right)q^n.\label{add*inter-proof-eqn=new}
\end{align}

By \eqref{proof-mk-0}, we derive that the coefficient of $q^n$ with $n\geq 0$ in \eqref{add*inter-proof-eqn=new} is nonnegative. The proof is complete.  \qed

\section{Proof of Theorem \ref{main-thm-proof-FHG}}

The objective of this section is to prove of Theorem \ref{main-thm-proof-FHG}. Throughout this section, we fix $t\geq 2$. For $0\leq i\leq t-1$,  define
\[F_i=q^{iR(t+1)-S(t-i)}(1-q^{R(t+1)-S})(1-q^{R(t+1)+S})(1-q^{S(2t+1-2i)}),\]
and for $-1\leq i\leq t-1$,  define
\[G_i=q^{(i+2)R(t+1)-S(t-i)}(1-q^{S(2t+1-2i)}).\]
We find that in order to prove Theorem \ref{main-thm-proof-FHG}, it suffices to show that
\begin{equation}\label{add-lem-3-1}
q^{-St} (1 - q^{R(t+1)-S})-q^{S(t+1)} (1 - q^{R(t+1)+S})=q^{Rt(t+1)}(1-q^S)+\sum_{i=0}^{t-1}F_i-G_{t-1},
\end{equation}
and for $0\leq i\leq t-1$,
\begin{equation}\label{add-lem-3-2}
F_i-G_i=H_i-G_{i-1}.
\end{equation}

We first give a proof of \eqref{add-lem-3-1}.

{\noindent \bf Proof of \eqref{add-lem-3-1}.} By elementary manipulations, we obtain
\begin{align*}
&\quad q^{-St} \left( 1 - q^{R(t+1)-S} \right) - q^{S(t+1)} \left( 1 - q^{R(t+1)+S} \right) \\
&= q^{-St} \left( 1 - q^{R(t+1)-S} \right) \left( 1 - q^{(t+1) \left(R(t+1)+S \right)} + q^{(t+1) \left(R(t+1)+S \right)} \right) \\
&\quad -q^{S(t+1)} \left( 1 - q^{R(t+1)+S} \right) \left( 1 - q^{(t+1) \left(R(t+1)-S \right)}+ q^{(t+1) \left(R(t+1)-S \right)} \right) \\
&= q^{-St} \left( 1 - q^{R(t+1)-S} \right) \left( 1 - q^{R(t+1)+S} \right) \sum_{i=0}^{t} q^{i \left(R(t+1)+S \right)} + q^{R(t+1)(t+1) + S} \left( 1 - q^{R(t+1)-S} \right) \\
&\quad- q^{S(t+1)} \left( 1 - q^{R(t+1)+S} \right) \left( 1 - q^{R(t+1)-S} \right) \sum_{i=0}^{t} q^{i \left(R(t+1)-S \right)}- q^{R(t+1)(t+1)} \left( 1 - q^{R(t+1)+S} \right)\\
&=\sum_{i=0}^{t} q^{iR(t+1)-S(t-i)}(1-q^{R(t+1)-S})(1-q^{R(t+1)+S})(1-q^{S(2t+1-2i)})\\
&\quad+q^{R(t+1)(t+1)}(q^S-q^{R(t+1)}-1+q^{R(t+1)+S})\\
&=\sum_{i=0}^{t-1}F_i +q^{Rt(t+1)}(1-q^{R(t+1)-S})(1-q^{R(t+1)+S})(1-q^S)-q^{R(t+1)(t+1)}(1-q^S)(1+q^{R(t+1)})\\
&=\sum_{i=0}^{t-1}F_i+q^{Rt(t+1)}(1-q^S)(1-q^{R(t+1)-S}-q^{R(t+1)+S}+q^{2R(t+1)}-q^{R(t+1)}-q^{2R(t+1)})\\
&=\sum_{i=0}^{t-1}F_i+q^{Rt(t+1)}(1-q^S)(1-q^{R(t+1)-S}-q^{R(t+1)+S}-q^{R(t+1)})\\
&=\sum_{i=0}^{t-1}F_i+q^{Rt(t+1)}(1-q^{R(t+1)-S}-q^{R(t+1)+S}-q^{R(t+1)}\\
&\qquad\qquad\qquad\qquad-q^S+q^{R(t+1)}+q^{R(t+1)+2S}-q^{R(t+1)+S})\\
&=\sum_{i=0}^{t-1}F_i+q^{Rt(t+1)}(1-q^{R(t+1)-S}-q^{S}+q^{R(t+1)+2S})\\
&=\sum_{i=0}^{t-1}F_i+q^{Rt(t+1)}(1-q^S)-q^{R(t+1)(t+1)-S}(1-q^{3S})\\
&=\sum_{i=0}^{t-1}F_i+q^{Rt(t+1)}(1-q^S)-G_{t-1}.
\end{align*}
This completes the proof.  \qed

Then, we give a proof of \eqref{add-lem-3-2}.

{\noindent \bf Proof of \eqref{add-lem-3-2}.} For $0\leq i\leq t-1$, we have
\begin{align*}
&\quad F_i-G_i\\
&=q^{iR(t+1)-S(t-i)}(1-q^{R(t+1)-S})(1-q^{R(t+1)+S})(1-q^{S(2t+1-2i)})\\
&\quad-q^{(i+2)R(t+1)-S(t-i)}(1-q^{S(2t+1-2i)})\\
&=q^{iR(t+1)-S(t-i)}(1-q^{S(2t+1-2i)})\left((1-q^{R(t+1)-S})(1-q^{R(t+1)+S})-q^{2R(t+1)}\right)\\
&=q^{iR(t+1)-S(t-i)}(1-q^{S(2t+1-2i)})(1-q^{R(t+1)-S}-q^{R(t+1)+S})\\
&=q^{iR(t+1)-S(t-i)}\\
&\quad\times(1-q^{R(t+1)-S}-q^{R(t+1)+S}-q^{S(2t+1-2i)}+q^{R(t+1)+S(2t-2i)}+q^{R(t+1)+S(2t+2-2i)})\\
&=q^{iR(t+1)-S(t-i)}(1-q^{R(t+1)+S}-q^{S(2t+1-2i)}+q^{R(t+1)+S(2t-2i)})\\
&\quad-q^{(i+1)R(t+1)-S(t-i+1)}(1-q^{S(2t+3-2i)})\\
&=H_i-G_{i-1}.
\end{align*}
The proof is complete.  \qed


\begin{thebibliography}{99}

\bibitem{Andrews-1976} G.E. Andrews, The Theory of Partitions, Encyclopedia of Mathematics and Its Applications, Vol. 2, Addison-Wesley, 1976.

\bibitem{Andrews-2019} G.E. Andrews and D. Newman, Partitions and the minimal excludant, Ann. Comb. 23(2019), 249--254.


\bibitem{Andrews-Merca-2012} G.E. Andrews and M. Merca, The truncated pentagonal number theorem, J. Comb. Theory Ser. A, 119 (2012) 1639--1643.

\bibitem{Andrews-Newman-2020} G.E. Andrews and D. Newman, The minimal excludant in integer partitions, J. Integer Seq. 23 (2020) Article 20.2.3.

\bibitem{Ballantine-Feigon-2012} C. Ballantine and B. Feigon, Truncated theta series related to the Jacobi Triple Product identity, Discrete Math. 348 (2025) 114319.

\bibitem{Ding-Sun-2025} X. Ding and L.H. Sun, Truncated theta series from the Bailey lattice, Adv. in Appl. Math. 167 (2025) 102884.


\bibitem{Ding-Sun-2025a} X. Ding and L.H. Sun, Proof of Merca's stronger conjecture on truncated Jacobi triple product series, (2025) arXiv: 2411.13818v3.

\bibitem{Grabner-Knopfmacher-2006} P.J. Grabner and A. Knopfmacher, Analysis of some new partition statistics, Ramanujan J. 12(3) (2006) 439--454.


\bibitem{Guo-Zeng-2013} V.J.W. Guo and J. Zeng, Two truncated identities of Gauss, J. Comb. Theory Ser. A, 120 (2013) 700--707.

\bibitem{Li-2023} X.L. Li, A generalized truncated sums of Jacobi triple product series and some related truncated theorems, Rev. Real Acad. Cienc. Exactas Fis. Nat. Ser. A-Mat. 117 (2023) 86.

\bibitem{Mao-2015} R. Mao, Proofs of two conjectures on truncated series,  J. Combin. Theory Ser. A 130 (2015) 15--25.



 \bibitem{Merca-2021} M. Merca, Rank partition functions and truncated theta identities. Appl. Anal. Discrete Math. (2021), DOI:10.2298/AADM190401023M.

 \bibitem{Merca-2021a} M. Merca, Truncated theta series and Rogers-Ramanujan functions, Exp. Math. 30(3) (2021) 364--371.

\bibitem{Wang-Yee-2019} C. Wang and A.J. Yee, Truncated Jacobi triple product series, J. Combin. Theory Ser. A 166 (2019) 382--392.

\bibitem{Xia-Zhao-2023} E.X.W. Xia and X. Zhao, Truncated sums for the partition function and a problem of Merca. Rev. Real Acad. Cienc. Exactas Fis. Nat. Ser. A-Mat. 116(1) (2022) 1--8.


\bibitem{Yee-2015} A.J. Yee, A truncated Jacobi triple product theorem, J. Combin. Theory Ser. A 130 (2015) 1--14.

\end{thebibliography}
\end{document}